\documentclass[10]{article}

\def\C{\mathbb{C}} 
\def\G{\mathbb{G}} 

\def\R{\mathbb{R}}  

\def\no{\noindent}

\def\ol{\overline}

\def\beq{\begin{equation}}
\def\eeq{\end{equation}}
\def\con{\overline}

\def\w{\wedge}

\usepackage{array}
\usepackage{tabularx}
\usepackage{amsfonts}
\usepackage{amssymb}
\usepackage{tikz-cd}
\usetikzlibrary{matrix,arrows,decorations.pathmorphing}
\usepackage{graphicx}   
\def\bpm{\begin{pmatrix}}
\def\epm{\end{pmatrix}}
\makeindex

\title{Geometric Matrices and the Symmetric Group}
\author{Garret Sobczyk \\
Universidad de las Am\'ericas-Puebla \\
 Departamento de F\'isico-Matem\'aticas \\
72820 Puebla, Pue., M\'exico}
	
\begin{document}

\maketitle

\begin{abstract} We construct $2^n \times 2^n$ real and complex matrices in terms of Kronecker products of a Witt basis of $2n$ null vectors in the geometric algebra $\G_{n,n}$ over the real and complex numbers. In this basis, every matrix is represented by a unique sum of products of null vectors. The complex matrices of $\G_{n,n+1}$  provide a direct matrix representation for geometric algebras $\G_{p,q}$, where $p+q\le 2n+1$. Properties of irreducible representations of the symmetric group are presented in this geometric setting. 
 
\smallskip
\no {\em AMS Subject Classification:} 15A66, 15B33, 81R05, 20C30
\smallskip

\no {\em Keywords: Clifford algebra, geometric algebra, geometric matrices, group algebra, representation theory, sparce matrices, symmetric group.}
\end{abstract}

 \section*{0\quad Introduction}
 
 A matrix is traditionally just a table of real or complex numbers. As
 eloquently put by Tobias Danzig, 
 \begin{quote} $\ldots$ a theory in which a whole
 	array of elements is regarded as a number-individual. These filing cabinets are added and multiplied, and a whole calculus of
 matrices has been established which may be regarded as a continuation of the algebra of complex numbers. \cite[p.212]{TD1967}. \end{quote}
 What is magical about matrices is contained in the definition of the {\it multiplication} of matrices. We give a geometrical construction of $2^n \times 2^n$ matrix algebras $Mat_{2^n}(F)$ over the real and complex numbers, when $F=\R$, or $F=\C$, respectively, in terms of the Kronecker Product of $2\times 2$-blocks of {\it anti-commuting} null vectors. 
 
 Whereas this rather restricted class of matrix algebras grow in size quite rapidly, 
  \[  Mat_{1}(F) , Mat_2(F), Mat_4(F) , \cdots, Mat_{2^n}(F) , \cdots \]
  they have unique geometric properties which call for special attention. Indeed,
  it is the geometric properties of the matrix algebras $Mat_2(\C )$ and $Mat_4(\C )$, that early in the $20^{th}$-Century was key to the development of the  fundamental theories of relativity and quantum mechanics. We explain in detail how the structure of the geometric algebras $\G_{p,q}$ are embedded in the matrix algebras $Mat_{2^n}(\C )$. The comprehensive geometric interpretation of matrices of these sizes has many ramifications, justifying the terminology {\it geometric matrices} \cite{S08}. One promising area of application is sparse random matrices, such as used in error correcting codes \cite{McKay}.  

The {\it permutation} on $n$ letters is a good place to start a study of the
symmetric group $S_n$, and the corresponding natural {\it permutation representation} of the rows of a square $n\times n$ matrix \cite{kao2010}. Such matrices are easily represented in terms of the geometric matrices of our approach. Representation theory of finite groups has provided new insight into the structure of finite groups, with many applications in other areas of mathematics and science. More generally, the study of finite dimensional representations of Lie groups and Lie algebra, \cite{JL2001,FH1991}, is a powerful tool in the study of quantum mechanics and atomic structure \cite{Weyl1950}. Aside from its many important applications, it is a beautiful mathematical structure that has developed over the last 150 years, and is worthy of study on its own merits.

\section{Geometric matrices}

Real or complex matrices of size $2^n\times 2^n$ are built up, block by block, by considering pairs of null vectors satisfying the following two simple rules:

  \begin{itemize}
  	\item[N1)] $a^2 = 0 = b^2$ \quad {\it The non-zero vectors $a $ and $b $ are nilpotents}.
  	\item[N2)] $a b + b a =1$ \quad {\it The sum of $a b $ and $b a $ is $1$.} 
   \end{itemize}
  Since we are assuming product associativity, the second property easily implies that $aba =a$.   
  We take $k$ such pairs of vectors $a_1, b_1, \cdots, a_k, b_k$, and assume that
  the null vectors with distinct indexes are pare-wise anti-commutative. That is,
  \[ a_i a_j = - a_j a_i , \ b_i b_j = - b_j b_i, \ a_i b_j = -  b_j a_i \]
  for $i \ne j$.
  
  If $M= \{m_{ij} \}$ and $N = \{ n_{jk} \}$ are matrices over geometric numbers, we define
  the {\it left} and {\it right} {\it Kronecker product} of $M$ and $N$ by
 \[ M \overrightarrow\otimes N:=  \{ M n_{ij}  \}, \ {\rm and} \ 
  M \overleftarrow\otimes N:=\{ m_{ij} N  \}.  \]
    As an example, consider the pairs of null vectors $a_i$ and $b_i$ satisfying N1) and N2) above. Then
  \[ \pmatrix{1 \cr b_1}\overrightarrow \otimes \pmatrix{1 \cr b_2} := 
  \pmatrix{1 \cr b_1 \cr b_2 \cr b_{12}}  \ {\rm and} \  \pmatrix{1 & a_2} \overleftarrow  \otimes \pmatrix{1 & a_1} := \pmatrix{1 & a_1 & a_2 & a_{21}}.    \]
  Defining $B_i^T:=\pmatrix{1 & b_i}$, and $A_i^T:= \pmatrix{1 & a_i}$, the above
  relationships take the form
  \[  B_1 \overrightarrow{\otimes}B_2=\pmatrix{B_1 \cr B_1 b_2} \ {\rm and} \ 
      A_2^T  \overleftarrow{\otimes}A_1^T = \pmatrix{A_1^T & a_2 A_1^T}.        \]
  Shortly, we will use these expressions and their generalizations.
  
  In a previous paper \cite{hypre2017}, I showed that the geometric
  algebra defined by
  \[ \G_{1,1}:=\Big\{  \R(a,b)| \ a^2 = b^2 = 0, \  ab+ba =1\Big\} , \] 
  is specified by its {\it spectral basis}
  \[ \G_{1,1} =B ab A^T=  \pmatrix{1 \cr b}ab \pmatrix{1 & a} =\pmatrix{ab  \cr b } \pmatrix{ab & a} = \pmatrix{ab & a \cr b & ba} .   \]
The matrix representation of a general element 
\[ g = g_{11}ab  + g_{12} a+g_{21} b + g_{22} ba \in \G_{1,1} \] 
is given by
\[  g = B^Tu A  g B^Tu A =       
     \pmatrix{1 & b}u\pmatrix{1 \cr a} g\pmatrix{1 & b}u \pmatrix{1 \cr a}            \] 
   \beq  = \pmatrix{1 & b}u \pmatrix{ g & g b \cr a g & agb} u \pmatrix{1 \cr a}=
      B^Tu [g] A ,    \label{specbasisg} \eeq
where $u=ab$ and $[g]:= \pmatrix{g_{11} & g_{12} \cr g_{21} & g_{22}}$ is the real matrix of the geometric number $g\in \G_{1,1}$ with respect to this spectral basis.
The matrix algebra $Mat_2(\R)$ is algebraically isomorphic, and thus fully equivalent to the geometric algebra $\G_{1,1}$, and inherits a comprehensive geometric interpretation.   

 By using the same spectral basis over the complex numbers, any geometric number $g\in \G_{1,2}$ can be
represented by a $2\times 2$ matrix $[g]_\C$ over the complex numbers. 
That is
\[  g =  \pmatrix{1 & b}ab [g]_\C  \pmatrix{1 \cr a} =B^Tu[g]_\C   A ,  \]
where the generating unit vectors   $e_1, f_1, f_2$ of $\G_{1,2}$ are defined by
\[ e_1=a+b, f_1 = a-b, \ {\rm and} \ f_2:= e_1 f_1 i, \]
and $i:=e_1 f_{12} $.

Calculating the conjugation operations \cite[p.60]{SNF}, of {\it reversion} $g^\dagger$ and {\it grade inversion} $g^-$, gives
\beq g^\dagger = A^Tu^\dagger \ol{[g]}^T  B    = B^T u\pmatrix{\ol g_{22} & \ol g_{12} \cr \ol g_{21} & \ol g_{11}}_\C  A 
  \label{dag-g}, \eeq
\beq    g^- =  \big( B^T\big)^- u\ol{[g]}  A^-      = B^Tu\pmatrix{\ol g_{11} & -\ol g_{12} \cr -\ol g_{21} & \ol g_{22}}_\C   A,  
\label{invertg} \eeq
respectively, 
and the reverse-inversion $g^*$ of $g \in \G_{1,2} $ is
\beq  g^* = \big(A^T\big)^-u^\dagger {[g]}^T  B^-   = B^T u\pmatrix{ g_{22} &- g_{12} \cr - g_{21} &  g_{11}}_\C  A.    \label{star-g} \eeq
 Using the operation of reverse-inversion, the {\it determinant} 
\beq \det [g] := gg^* = B^Tu [g][g^*] A = g_{11} g_{22}- g_{12} g_{21} . \label{detg} \eeq  

Before generalizing to higher dimensional geometric algebras $\G_{p,q}$, consider $\G_{2,2}$. The spectral
basis of $\G_{2,2}$ is  
\[  \G_{2,2} =B_1 \overrightarrow{\otimes}B_2 u_{12} A_2^T \overleftarrow{\otimes}A_1^T =\pmatrix{B_1 \cr B_1 b_2}u_{12} \pmatrix{A_1^T& a_2 A_1^T}\] 
\[=
\pmatrix{B_1 u_{12}A_1^T & B_1 u_{1}a_2 A_1^T \cr B_1 b_2 u_{1}A_1^T &
	B_1 u_{1}u_2^\dagger  A_1^T }\]
\[=\pmatrix{1 \cr b_1}\overrightarrow\otimes \pmatrix{1 \cr b_2}
u_{12}  \pmatrix{1 & a_2} \overleftarrow \otimes \pmatrix{1 & a_1 }  \]  
 \[ =\pmatrix{1 \cr b_1}\overrightarrow\otimes \pmatrix{u_{12} & a_2 u_1 \cr b_2 u_1 & u_1 u_2^\dagger} \overleftarrow \otimes  \pmatrix{1 & a_1}  \] 
 \beq = \pmatrix{u_{12} & a_{1} u_2 & a_{2}u_1 & a_{21} 
 		\cr b_1 u_2 & u_1^\dagger u_2 & b_1 a_2 &-a_2 u_1^\dagger   \cr
 		b_{2}u_1  & b_2 a_1 & u_1 u_2^\dagger & a_1 u_2^\dagger  \cr
 		b_{12} & -b_2 u_1^\dagger & b_1u_2^\dagger & u_{12}^\dagger} \label{specG22} \eeq  	    
 	where $u_i:= a_ib_i$, $u_i^\dagger=b_i a_i$, $u_{12}:=u_1u_2$, and $u_{12}^\dagger = u_1^\dagger u_2^\dagger$.

For $g\in \G_{2,2}=\G_{1,1}\times \G_{1,1}^\prime, $
where $\G_{1,1}:=\R(a_1,b_1)$ and $\G_{1,1}^\prime:=\R(a_2,b_2)$,  
   \[ g= B_1^T \overrightarrow{\otimes}B_2^Tu_{12} \Big[  A_2 \overleftarrow{\otimes}A_1 g
   B_1^T \overrightarrow{\otimes}B_2^T  \Big] u_{12}A_2 \overleftarrow{\otimes}A_1     \]  
   \[  = B_1^T \overrightarrow{\otimes}B_2^Tu_{12} [ g ] A_2 \overleftarrow{\otimes}A_1 = B_1^T u_1\Big[ B_2^Tu_2\pmatrix{[g]_{11} & [g]_{12} \cr [g]_{21} & [g]_{22}}A_2 \Big]A_1               \]
   \beq   =B_1^T u_1\pmatrix{g_1^\prime & g_2^\prime \cr g_3^\prime & g_4^\prime}A_1              \label{2blocks} \eeq
   where $[g]$ is the real $4\times 4$ matrix of $g$ with respect to the spectral basis (\ref{specG22}), $[g]_{kl}$ are its $2\times 2$ blocks, and $g_1^\prime,g_2^\prime, g_3^\prime, g_4^\prime$ are the elements in $\G_{1,1}^\prime$ represented by each of these blocks, respectively.  Generalizing (\ref{2blocks}) for $\G_{n+1,n+1}=\G_{1,1}\times \G_{n,n}^\prime $, gives for $g\in \G_{n+1,n+1}$  
      \beq g =B_1^T u_{1}\pmatrix{g_1^\prime & g_2^\prime \cr g_3^\prime & g_4^\prime}A_1,  \label{specGnn} \eeq  
   for $g_1^\prime,g_2^\prime,g_3^\prime,g_4^\prime \in \G_{n,n}^\prime$. 
   
    For $g\in \G_{2,2}=\G_{1,1}\times \G_{1,1}^\prime$, the reverse of $g$ is
   \[g^\dagger = A_1^T \overrightarrow{\otimes}A_2^Tu_{12}^\dagger [ g]^T
   B_2 \overleftarrow{\otimes}B_1 =\pmatrix{A_1^T & A_1^Ta_2}u_{12}^\dagger[g]^T \pmatrix{B_1 \cr b_2 B_1}.       \]  
   We also have 
   \[ [g] = \pmatrix{g_{11} & g_{12} & g_{13} & g_{14} \cr    
   	g_{21} & g_{22} & g_{23} & g_{24} \cr  
   	g_{31} & g_{32} & g_{33} & g_{34} \cr  
   	g_{41} & g_{42} & g_{43} & g_{44} }, \ \ 
   	 [g^\dagger] = \pmatrix{g_{44} & g_{34} & -g_{24} & -g_{14} \cr    
   	 	g_{43} & g_{33} & -g_{23} & -g_{13} \cr  
   	 	-g_{42} & -g_{32} & g_{22} & g_{12} \cr  
   	 	-g_{41} & -g_{31} & g_{21} & g_{11} }, \ \ 
   	   	\]
   	 and
   \[  [g^*] = \pmatrix{g_{44} & -g_{34} & g_{24} & -g_{14} \cr    
   -	g_{43} & g_{33} & -g_{23} & g_{13} \cr  
   	g_{42} & -g_{32} & g_{22} & -g_{12} \cr  
   	-g_{41} & g_{31} & -g_{21} & g_{11} }, \ \ 
   \]

 The matrices $[a_1],[b_1],[a_2],[b_2]$ of the null vectors $a_1,b_1,a_2,b_2$, are

 \[  [a_1]= \pmatrix{0 & 1 & 0 & 0 \cr
 	                 0 & 0 & 0 & 0 \cr
 	                 0 & 0 & 0 & 1 \cr 
 	                 0&0&0&0 }, \ 
 	   [a_2] =     \pmatrix{0 & 0 & 1 & 0 \cr
 	                 0 & 0 & 0 & -1 \cr
 	                 0 & 0 & 0 & 0 \cr 
 	                 0&0&0&0 }   \]	
  \[  [b_1]= \pmatrix{0 & 0 & 0 & 0 \cr
  	1 & 0 & 0 & 0 \cr
  	0 & 0 & 0 & 0 \cr 
  	0&0&1&0 }, \ 
  [b_2] =     \pmatrix{0 & 0 & 0 & 0 \cr
  	0 & 0 & 0 & 0 \cr
  	1 & 0 & 0 & 0 \cr 
  	0&-1&0&0 }  .  \]
  Note that in the special case when 
  \[ g=g_{11} a_1b_1+g_{12} a_1+g_{21}b_1+g_{22}b_1 a_1\in \G_{1,1}\subset \G_{1,1}\times \G^\prime_{1,1}  , \]
   then the $4\times 4$-matrix $[g]$ of $g$ is
  \[  [g]=\pmatrix{g_{11}&g_{12} & 0 &0 \cr g_{21} & g_{22} & 0& 0 \cr 0&0 & g_{11}&g_{12}  \cr 0&0& g_{21} & g_{22} } .  \]
  Alternatively, for $g^\prime =g_{11}^\prime a_2 b_2 + g_{12}^\prime a_2+ g_{21}^\prime b_2+g_{22}^\prime b_2a_2\in \G_{1,1}^\prime$, then  	
    \[  [g]=\pmatrix{g_{11}^\prime & 0 &g_{12}^\prime & 0  \cr
    	 0 & g_{11}^\prime & 0& - g_{12}^\prime \cr
    	 g_{21}^\prime &0 & g_{22}^\prime & 0  \cr
         0 & -g_{21}^\prime & 0 & g_{22}^\prime } .  \]
    
    The {\it standard basis} of $\G_{n,n}$ is defined by
    \[ \G_{n,n} := \R(e_1, \cdots, e_n,f_1, \cdots f_n)
    =gen_\R \{e_i,f_i| \ e_i := a_i+b_i, f_i=a_i-b_i \},     \]
    for $1\le i \le n$.  
    For the {\it complexified geometric algebra} 
    \[ \G_{n,n+1}:= \C(e_1,\cdots,e_n,f_1,\cdots f_n)=gen_\C \{e_i,f_i| \ e_i := a_i+b_i, f_i=a_i-b_i \} ,\]
    with $i:= e_{n\cdots 1}f_{1\cdots n+1}$. Defined in this way, the $2n+1$-vector $i$, playing the roll of the imaginary unit, is in the center of $\G_{n,n+1}$, and $i^2=-1$.     
    
  Consider now the geometric algebra $\G_{p+1,q+1}=\G_{p,q}^\prime\times \G_{1,1}$. Using the spectral basis of $\G_{1,1}=gen_\R\{a_1,b_1\}$, a general element $g\in \G_{p+1,q+1}$ can be written
  \[  g= h_1 a_1b_1 + h_2 a_1 + h_3 b_1 + h_4 b_1a_1 = a_1b_1h_1 +  a_1 h_2^- +  b_1 h_3^- +  b_1a_1 h_4 , \]
  where $h_1,h_2,h_3,h_4 \in \G_{p,q}^\prime$. Then, for $u_1=a_1b_1$,
  \beq  g = \pmatrix{1 & b_1}u_1 \pmatrix{1 \cr a_1 }g \pmatrix{1 & b_1}u_1  \pmatrix{1 \cr a_1}      
  =  B^T u_1\pmatrix{h_1 & h_2 \cr h_3^- & h_4^- } A. \label{pqdecomp} \eeq
   Thus, the elements of $\G_{p+1,q+1}$ can be expressed as $2\times 2$ matrices over
  $\G_{p,q}^\prime$.

The spectral basis referred to in \cite[p.206]{SNF} is equivalent to
\[ \pmatrix{1 & a_1 & a_2 & a_{21}}^T b_1 a_1 b_2 a_2 \pmatrix{1 & b_1 & b_2 & b_{12}}   \]	
The {\it reverse $\dagger$} of this basis is
\[   \pmatrix{1 & b_1 & b_2 & b_{21}}^T a_1 b_1 a_2b_2 \pmatrix{1 & a_1 & a_2 & a_{12}} ,   \]
which is {\bf not} equivalent to the spectral basis (\ref{specG22}), meaning there is no inner automorphism which will take this spectral basis into the spectral basis (\ref{specG22}), and at the same time leave invariant of the null vectors
$a_i$ and $b_i$.

 \section{The permutation group algebra of $ S_n $}
  
    The symmetric group $S_n$ consists of all permutations on $n$ letters. We use the
  usual cycle notation. For example,
  \[ S_3 =\{ 1, (12),(13),(23),(123),(132)\},   \]
  with group multiplication of cycles. An example cycle multiplication, from right to left, is \[ (23)=(12)(13)(12) .\]

The {\it permutation representation} of $S_n$ is generated by real $n\times n$ matrices of the form
 \[  (1k)={\bf 1}_{(1k)}, \]
 where the matrix ${\bf 1}_{(1k)}$ is obtained by interchanging the first and the
 $k^{th}$ rows of the $n\times n$ of the identity matrix $\bf 1$. The generators of the permutation matrix group algebra ${ S}_3$ are
 \[{ S}_3 =gen \big\{(12),(13) \big\} = gen\Big\{ \pmatrix{0 & 1 & 0 \cr 1 & 0 & 0 \cr 0 & 0 & 1}, \pmatrix{0 & 0 & 1 \cr 0 & 1 & 0 \cr 1 & 0 & 0} \Big\}. \]
 
 In studying the permutation representation, and the closely related {\it standard representation} of $S_n$, it is expedient to introduce the {\it all ones matrix}
 \[ {\cal A}_n:= \pmatrix{1 & \cdots &1 \cr \cdot & \cdots & \cdot
 	 \cr \cdot & \cdots & \cdot \cr
 	            1 & \cdots & 1},   \]
 and the {\it Casimir matrix} \cite{Kun2015},
 \[  {\cal C}_n := {\cal A}_n - 1_n,  \] 
 where $1_n$ is the $n\times n$-identity matrix. For example
 \[ {\cal C}_3 = \pmatrix{0 & 1&1 \cr 1 & 0 & 1 \cr 1 & 1 & 0}.    \]
 
 Not surprisingly, both the $n$-all ones matrix and the $n$-Casimir matrix commutes with all the regular symmetric matrices $ S_n$. It is also easy to verify the simple relationships
 \[ {\cal A}_n^2 = n {\cal A}_n \ \ {\rm and} \ \ {\cal C}_n^2=  (n-2){\cal C}_n +(n-1)1_n ,  \]
 implying that
  \beq  min({\cal A}_{n}) := x(x-1) , \ {\rm and} \ min({\cal C}_{n}) :=(x+1)\big(x-(n-1)\big),  \label{minpolys}  \eeq
  are, respectively, the minimal polynomials of ${\cal A}_n$ and ${\cal C}_n$.
  	  
  \section{The geometric group algebra of $S_n$}            
 	                
The geometric algebras $\G_{n,n}$ and $\G_{n,n+1}$ provide a comprehensive geometric interpretation for the elements of the group algebra of the symmetric group $ S_n$, and give a new
way of studying its properties. One problem of this approach is that these geometric group algebras only exist for dimensions $2^n$, corresponding to the matrix
representations using $2^n \times 2^n$ matrices over the real or complex numbers. However, our approach provides both new tools and a geometric interpretation of the results.   

We begin by examining the geometric group algebras of the smaller symmetric groups, and then generalize these results to $S_n$. The group algebra of
$S_2$ is generated by the single element $(12)=\pmatrix{0&1 \cr 1 & 0}$. The
corresponding generating element in the geometric group algebra is $(12)=a_1 +b_1 \in \G_{1,1} $. Thus, the permutation representation of $S_2$ is
\[ S_2 = span \{ \pmatrix{1& 0 \cr 0 & 1} ,\pmatrix{0&1 \cr 1 & 0} \}=span_\R \{1, a_1+b_1\}=gen_{\R } \{ a_1+b_1 \}         . \] 

To get the permutation geometric group algebra of $S_3$, we must go to the larger geometric algebra $\G_{2,2}$. Using the spectral basis (\ref{specG22}), we find the geometric permutation representation
\beq S_3 =gen_\R \{(12),(13)\} = gen_\R \{1+(a_1+b_1 -1)u_2, 1+ (a_2+b_2-1)u_1 \}. \label{permS3a} \eeq
For $S_4$, we have
\beq S_4=S_3 \cup \{(14)\}=S_3 \cup \{ 1+ (a_{21}+ b_{12} - u_{12} - u_{12}^\dagger) \}.         \label{permS3b} \eeq 
In a different approach, \cite[p.201-222]{SNF}, a special {\it twisted product} was developed, but has not yet been fully explored.

There is an {\it irreducible representation}, called the {\it standard representation}, which is easily obtained from the permutation representation of $S_n$. We illustrate the general method first for $S_5$ in $\G_{2,2}$.
The minimal polynomial (\ref{minpolys}) for ${\cal C}_4$ is  
\[ min({\cal C}_4)= (x+1)(x-3).  \]
The {\it spectral basis} \cite[p.125]{SNF}, \cite{S2,S0}, for this minimal polynomial is 
\[ s_1=\frac{{\cal C}_4-3}{-4}=\frac{{\cal A}_4-4}{-4}, \ \ {\rm and} \ \ s_2=\frac{{\cal C}_4+1}{4}=\frac{{\cal A}_4}{4},   \]
where $s_1$ and $s_2$ are mutually annihilating idempotents with the property that $s_1+s_2=1$.

Using (\ref{specG22}), from the matrices ${\cal A}_4$ and ${\cal C}_4$, we calculate
the corresponding geometric numbers $A_4$ and $C_4$ in $\G_{2,2}$, getting
\[ A_4 = 1+a_1+b_1 +(a_2+b_2)\Big( (a_1-b_1)+2 a_1 \w b_1 \Big)    \] 
\[ = A_2\Big( 1+2(a_2+b_2)\Big)a_1\w b_1 =1+ C_4.  \]
where $A_2:= (1+a_1+b_1) \in \G_{1,1}\subset \G_{2,2}$.

More generally, for the minimal polynomials (\ref{minpolys}), the spectral basis
of $C_{2^n}\in \G_{n,n} $ is given by
 \beq s_1=\frac{C_{2^n}-2^n+1}{-2^n}=\frac{A_{2^n}-2^n}{-2^n}, \ \ {\rm and} \ \ s_2=\frac{C_{2^n}+1}{2^n}=\frac{A_{2^n}}{2^n}.   \label{specCn} \eeq
  We have the useful 
 recursive relation
 \[  A_{2^{n+1}} = A_{2^n}\Big(1+2^n(a_{n+1}+b_{n+1})(a_1\w b_1) \cdots (a_{n}\w b_{n})\Big) =1+C_{2^{n+1}} \in \G_{n+1,n+1},\]
 where $A_{2^n} \in \G_{n,n}$ for $n\ge 1$. For $n=1$, $A_1 =1$ and $C_1=0$.
 
 Using the spectral basis (\ref{specCn}),
 \[  C_{2^n} = (-1) s_1 + (2^n-1)s_2. \]
Since the Casimir geometric number $C_{2^n}$ commutes with the permutation representations of $S_{2^n}$, 
it follows that the inner automorphism of the geometric number 
\beq   g_c=s_1(1-u_{1\cdots n}^\dagger) +s_2u_{1\cdots n}^\dagger \label{gcsurgery} \eeq
 will {\it diagonalize} $C_n$, meaning that 
 \[ [g_c^{-1}C_n g_c ] =\pmatrix{-1 & 0 &\cdots & 0 &0 \cr 
 	        0 & -1 &\cdots & 0 &0\cr  
 	         &\cdots  & \cdots &  \cr
 	          & \cdots &\cdots &  \cr
 	        0 & 0 &\cdots &- 1 & 0 \cr 
 	        0 & 0 &\cdots &  0 & 2^n-1             }.    \]
 If we now apply the inner automorphism of $g_c$ to the permutation representation
 of $S_{2^n}$ in $\G_{n,n}$, we get the standard irreducible representation of
 $S_{2^{n}} $ in $\G_{n,n}$.
 
 The explicit calculations for the standard irreducible representation of $S_4$ in $\G_{2,2}$ are as follows:
  \[ C_4 = (-1)s_1 + 3 s_2 \] 
  for $s_1 = \frac{C_4 -3}{-4}$, and $s_2= \frac{C_4+1}{4}$. Next, we preform surgery,
  removing the last column of $s_1$, and replacing it with the last column of $s_2$, giving 
  \[ g_c = s_1(1-u_{12}^\dagger)+ s_2 u_{12}^\dagger ,  \]
which is a matrix whose column vectors are the eigenvectors of $C_4$. Checking,
  \beq [g_c^{-1} C_4 g_c] = \pmatrix{-1 & 0 & 0 & 0 \cr
  	                               0 & -1 & 0 & 0 \cr
  	                               0 & 0 & -1 & 0 \cr
  	                               0 & 0 & 0 & 3 }, \label{gcsurgery2}  \eeq
 showing the $g_c$ diagonalizes $C_4$ as expected. 	                            
 
 Applied to the permutation representations $(12),(13)$ and $(14)$ of $S_4$ in $\G_{2,2}$, given in (\ref{permS3a}) and (\ref{permS3b}), $g_c$ {\it block
 diagonalizes} the matrix generators of these permutations, giving the {\it standard irreducible representation} of $S_4$,
 \[  (12)=\pmatrix{0 & 1 & 0 & 0 \cr
 	1 & 0 & 0 & 0 \cr
 	0 & 0 & 1 & 0 \cr
 	0 & 0 & 0 & 1 } =[1+(a_1+b_1-1)u_2]   ], \] 
  \[  (13)=\pmatrix{0 & 0 & 1 & 0 \cr
                          	0 & 1 & 0 & 0 \cr
                        	1 & 0 & 0 & 0 \cr
                          	0 & 0 & 0 & 1 } =[1+(a_2+b_2-1)u_1]\]
 and 
 \[  (14)=\pmatrix{-1 & 0 & 0 & 0 \cr
             	    -1 & 1 & 0 & 0 \cr
 	                -1 & 0 & 1 & 0 \cr
 	                0 & 0 & 0 & 1 }=[1-(2+b_1+b_2)u_{12}].     \]
Note, in this representation the permutation matrices of $(12)$ and $(13)$ are unchanged, but the representation of $(14)$ is now a $3\times 3$ sub-block.  	              	  

From what we have learned, it is easy to write down a standard representation of $S_5$, completing the permutation representations for the $2$-cycles $(12),(13,(14)$, 
with
 	          \[   (15) = \pmatrix{-1 & 0 & 0 & 0 \cr
 	          	-1 & 1 & 0 & 0 \cr
 	          	-1 & 0 & 1 & 0 \cr
 	          	-1 & 0 & 0 & 1   }=[1-(2+b_1+b_2+b_{12})u_{12}] .    \]  
 	          
\section{Group characters}

In studying a particular group algebra, defined by its permutation representation as a subgroup of the symmetric group, it is very useful to find a matrix which commutes 
with all the conjugacy classes of the symmetric group, because these elements are one way of finding the respective irreducible representations \cite{newman67}. Here we study the conjugacy classes of the group algebra of $S_4$.  

There are five conjugacy classes of $S_4$:
The group identity $1$, the six $2$-cycles, the $8$ 3-cycles, the 6 $4$-cycles, and 3 double $2$-cycles, making up the 24 group elements of $S_4$. The general matrix that commutes with the generating 2-cycles $(12),(13)$, and $(14)$, of $S_4$ will necessarily commute with all of the elements of $S_4$. The matrix 
$[g_{ij}]$ that commutes with $S_4$ is 
\beq g_{all}:=\pmatrix{s & t & t& t \cr   
	          t & s & t& t \cr
	          t & t & s& t \cr
	          t & t & t& s },     \label{22D3CY} \eeq
dependent on the two independent parameters $s,t$. The chacteristic and minimal polynomials of $g_{all}$ are $\Big(\lambda-(s-t)\Big)^3\Big( \lambda-(3t+s)\Big) $ and
\[ min(g_{all})=\Big(\lambda-(s-t)\Big)\Big( \lambda-(3t+s)\Big) , \]
respectively.

 The matrix $g_{alt}$ which commutes with the conjugacy class of the alternating group, or even subgroup of $S_4$ is  
 \beq g_{alt}=\pmatrix{s_1 & t_1 & s_2& t_2 \cr   
 	t_1 & s_1 & t_2& s_2 \cr
 	s_2 & t_2 & s_1& t_1 \cr
 	t_2 & s_2 & t_1& s_1 },   \label{4CY} \eeq
dependent on the four parameters $s_1, t_1, s_2,t_2$. The minimal polynomial of
$g_{alt}$ is
\[  min(g_{alt})=\Big(\lambda -(t_2-s_2-t_1+s_1)\Big)\Big(\lambda- (-t_2+s_2-t_1+s_1)\Big)\] \[\Big(\lambda-(-t_2-s_2+t_1+s_1)\Big) \Big(\lambda- (t_2+s_2+t_1+s_1)\Big).    \] 	            
 Clearly, the condition (\ref{4CY}) of commuting with all even-cycles is less restrictive  than the condition (\ref{22D3CY}) of commuting with all of $S_4$.
 
 Representation theory is built up on an {\it inner product} defined on the {\it traces} of matrices which define the {\it characters} of a representation. Given an element
 $g \in \G_{n,n+1}$, together with its formally {\it complex matrix} $[g]$, 
 \[  tr[g] := 2^n \langle g \rangle_{0+(2n+1)}, \quad {\rm and} \quad
      tr\con{[g]} := 2^n \langle g^- \rangle_{0+(2n+1)},           \] 
 where $i:= e_{1\cdots n}f_{1 \cdots n+1}^\dagger$. A {\it representation} of a group $G$ in
 the geometric algebra $\G_{n,n+1}$ is a mapping 
 \[    f: G \to \G_{n,n+1}\] 
  of the group $G$ into the geometric algebra $\G_{n,n+1}$ which preserves 
  the group operation, and for which $f_e = 1$, where $e \in G$ is the group identity. It follows that $f_{gh} = f_g f_h$ for all $g,h \in G$.   
  Let $f$ be a representation of the group $G$ in $\G_{n,n+1}$. The
  {\it group character} 
  \[   \psi_f: G \to \C  \]
  is defined by $ \psi_f(g):=tr(f_g)$ for each $g\in G$. A comprehensive introduction to representation theory is found in \cite{FH1991}.

 \section{New tools}	                
 	            
 Consider an arbitrary geometric number $g\in \G_{2,2}$, with the matrix
 $[g_{ij}]$ in the basis (\ref{specG22}). We can perform ``surgery" on this matrix. For example,
    \[  [g-g u_{2}^\dagger-u_2^\dagger g] = \pmatrix{g_{11} & g_{12} & 0 & 0 \cr
    	                                            g_{21} & g_{22} &0 & 0 \cr
    	                                            0 & 0 &-g_{33} & 0 \cr 
    	                                            0 &0&0&g_{44}}     \]	 
 and 
  \[ [g-g u_{12}^\dagger-u_{12}^\dagger g] = \pmatrix{g_{11} & g_{12} & g_{13} & 0 \cr
   	g_{21} & g_{22} &g_{23} & 0 \cr
   	g_{31} & g_{32} &g_{33} & 0 \cr 
   	0 &0&0&-g_{44}}.     \]
      The same such surgery was performed on $s_1$ and $s_2$ in (\ref{gcsurgery}) and
   (\ref{gcsurgery2}) to obtain the respective automorphisms defined by $g_c$.
   We can also directly extract rows or columns of the matrix $[g]$. For example,
   \[ [g b_{1}u_2]= \pmatrix{g_{12} &0&0&0 \cr
   	                      g_{22} &0&0&0 \cr
   	                      g_{32} &0&0&0 \cr
   	                      g_{42}&0&0&0} \quad {\rm and} \quad 
   	          [gb_{12}]= \pmatrix{g_{14} &0&0&0 \cr
   	          	g_{24} &0&0&0 \cr
   	          	g_{34} &0&0&0 \cr
   	          	g_{44} &0&0&0}           \]  
   
   Direct representations of matrix algebras by Clifford algebras $\G_{n,n}$  of neutral signature restricts the classes of matrices considered to sizes of $2^n\times 2^n$. However, the operation of surgery ameliorates this restriction by allowing us to extract matrices of smaller sizes. Because of their geometric significance, matrices over the complex numbers of size $2\times 2$, the {\it Pauli Matrices}, and of size $2^2\times 2^2$, the {\it Dirac Matrices}, have proven themselves to be of immense importance in physics. It is this extra geometric structure that geometric algebra makes explicit in its many different guises \cite{geoS2017}. 	           
 	        
 The {\it standard spectral basis} of the geometric algebras $\G_{n,n}$ and $\G_{n,n+1}:=\G_{n,n}(\C)$ are easily constructed, starting with the primitive idempotent $u_{1\cdots n}:=u_1\cdots u_n$
 where each $u_i = a_i b_i$. The first column is then written down, followed by successive columns defined by the Kronecker products, with the dual blocks in reverse order, 
 \[ \G_{n,n}:= B_1 \overrightarrow{\otimes} \cdots \overrightarrow{\otimes} B_n u_{1\cdots n} A_n^T \overleftarrow{\otimes} \cdots \overleftarrow{\otimes} A_1^T. \] 
 We have already demonstrated the method for $n=2$ in (\ref{specG22}); we now
 give the standard spectral basis for $\G_{3,3}$.  
 
 The standard spectral basis defining $\G_{3,3}$ is
 \[ \G_{3,3}:= B_1 \overrightarrow{\otimes} B_2 \overrightarrow{\otimes} B_3 u_{123}
  A_3 \overleftarrow{\otimes} A_2 \overleftarrow{\otimes} A_1.    \]
  In expanded form, $\G_{3,3}=$ 
 \beq \pmatrix{u_{123} & a_1 u_{23} & a_2 u_{13} & a_{21} u_{3} & a_3 u_{12}
 	            & a_{31} u_2 & a_{32} u_1 & a_{321}  \cr
 	            b_1 u_{23} & u_1^\dagger u_{23} & b_1 a_2 u_{3} & -a_{2}u_1^\dagger u_{3} & b_1 a_3 u_{2}
 	            & -a_{3} u_1^\dagger u_2 & b_1 a_{32} & a_{32}u_1^\dagger \cr
 	          b_{2} u_{13} & b_2 a_1 u_{3} & u_2^\dagger u_{13} & a_{1}u_2^\dagger u_{3} & b_2 a_3 u_{1}^\dagger
 	            & b_2 a_{31} & b_2 a_{32} u_1 & a_{13} u_{2}^\dagger \cr
 	           b_{12} u_{3} & -b_2 u_1^\dagger u_{3} & b_1 u_2^\dagger u_{3} & u_{12}^\dagger u_{3} & a_3 b_{12}
 	            &b_2 a_{3} u_1^\dagger & a_{3}b_1  u_2^\dagger & a_{3}u_{12}^\dagger \cr
 	            b_3 u_{12} & b_3 a_1 u_{2} & b_3 a_2 u_{1} & b_3 a_{21} & u_{12} u_{3}^\dagger
 	            & a_{1} u_2 u_3^\dagger & a_2 u_{1} u_3^\dagger & a_{21} u_3^\dagger  \cr
 	           b_{13} u_{2} & -b_3 u_1^\dagger u_{2} &b_{13} a_2 & b_3 a_{2} u_{1}^\dagger & b_1 u_{2} u_3^\dagger 
 	            & u_{13}^\dagger u_2 & b_1 a_{2} u_3^\dagger & -a_2 u_{13}^\dagger \cr
 	           b_{23} u_{1} & b_{23} a_1 & -b_3 u_1 u_{2}^\dagger & -b_3 a_{1} u_{2}^\dagger & b_2 u_1 u_{3}^\dagger
 	            & b_2 a_{1} u_3^\dagger & u_1 u_{23}^\dagger & a_1 u_{23}^\dagger \cr
 	            b_{123} & b_{23} u_{1}^\dagger & -b_{13} u_{2}^\dagger & b_{3} u_{12}^\dagger & b_{12} u_{3}^\dagger 
 	            & -b_2 u_{13}^\dagger & b_1 u_{23}^\dagger & u_{123}^\dagger  }   \label{stdspecG33} \eeq
 
  We can now easily construct the geometric permutation representation of
  $S_9$ in $\G_{3,3}$. Several examples of the geometric permutations are
  \[ (12) = 1+ (a_1+b_1-1)u_{23}, \ \ (13)=1+(a_2+b_2-1)u_{13} ,  \] 
   \[ (16)= 1+ (a_{31}+b_{13} -u_{13}-u_{13}^\dagger)u_2   \]
   and
 \[ (19)= 1-u_{123} - (1 +b_1)(1+b_2)(1+b_3)u_{123}. \] 
  It is interesting to write down the matrix for the 9-cycle $(123456789)$. We have
 \[   (123456789) = \pmatrix{0& 0&0&0&0&0&0&-1 \cr   
                            1& 0&0&0&0&0&0&-1 \cr
                            0& 1&0&0&0&0&0&-1 \cr
                            0& 0&1&0&0&0&0&-1 \cr
                            0& 0&0&1&0&0&0&-1 \cr
                            0& 0&0&0&1&0&0&-1 \cr
                            0& 0&0&0&0&1&0&-1 \cr
                            0& 0&0&0&0&0&1&-1 }    \] 
      \[=[b_1+b_2a_1+b_3 a_{21} - (a_{321} +a_{32}+a_3 +a_{21}-a_2+a_1+1)u_{123}^\dagger],          \]
  which is a $9^{th}$ root of unity.  
  
  Let us use our new tools to decompose the regular representation of $S_3$ in $\G_{3,3}$ into the sum of its irreducible parts \cite[p.127]{JL2001}. We begin by writing
  \[ X=x_0 + x_1 (18)+x_2 (19) +x_3 (89) + x_4 (189) + x_5 (198), \]                       
  represented by an $8\times 8$ matrix $[X]$. We have used the generating $2$-cycles $(18)$ and $(19)$ to take advantage of the irreducible standard representation of $S_3$ in $Mat_2(\R)$. Letting $x_{0\cdot 5}=x_0 + \cdots + x_5$, the matrix $[X]$ of $X$ is
  \[ \pmatrix{x_0-x_2+x_3-x_5& 0& 0 &0&0&0&0&x_1-x_3-x_4+x_5 \cr   
                -x_2-x_5&x_{0\cdot 5}& 0 &0&0&0&0&0 \cr 
               -x_2-x_5& 0 & x_{0\cdot 5} &0&0&0&0&0 \cr 
                -x_2-x_5& 0&0&x_{0\cdot 5}&0&0&0& 0 \cr
                -x_2-x_5& 0&0&0&x_{0\cdot 5}&0&0& 0 \cr
                -x_2-x_5& 0&0&0&0&x_{0\cdot 5}&0& 0 \cr
                -x_2-x_5& 0&0&0&0&0&x_{0\cdot 5}& 0 \cr
                	x_1-x_2+x_4-x_5& 0&0&0&0&0&0& x_0+x_2-x_3-x_4 }.              \]
  
   Except for the strange arrangement of $2\times 2$ block, the block diagonal matrix $[X]$ has desired form. The matrix $[X]$ can be fully diagonalized, however this would not give the desired irreducible decomposition of $[X]$. Instead, noting that $x_{0\cdot 5}$ is an eigenvalue of the remaining $2$-block, we block diagonalize $X$ with the eigenvectors of the matrix $[X]_{x_0\to 1, x_1 \to 2, \cdots, x_5 \to 6}$,
    getting the decomposition of $[X]$ into the sum of its irreducible representations
   \[  \pmatrix{x_{0\cdot 5}& 0&0&0&0&0&0&0 \cr   
   	 0& x_{0\cdot 5}&0 &0&0&0&0&0 \cr 
   	 0& 0 &x_{0\cdot 5}&0&0&0&0&0 \cr 
   	0& 0&0&x_{0\cdot 5}&0&0&0& 0 \cr
   	0& 0&0&0&x_{0\cdot 5}&0&0& 0 \cr
   	0& 0&0&0&0&x_{0\cdot 5}&0& 0 \cr
   	0& 0&0&0&0&0&	x_0+x_1-x_3-x_4& x_2-x_3-x_4+x_5 \cr
   	0& 0&0&0&0&0&-x_1+x_2+x_4-x_5& 	x_0-x_1+x_3-x_5 }.          \] 
   Since we are representing the group algebra of $S_6$ in $\G_{3,3}$, we can
   reduce the size of our representation by acting on the $6$-dimensional column vector
   \[ x=\pmatrix{0 & 0 & x_0 & x_1 & x_2 & x_3 & x_4 & x_5}^T, \]
    which essentially eliminates the first two rows and columns of the matrix.
    
    \section*{Appendix: Lower dimensional geometric algebras} 
    
    Recall that the standard basis of the geometric algebras $\G_{n,n}$ and
    $\G_{n,n+1}$ are defined by
    \[ \G_{n,n}:= \R(e_1, \cdots , e_n,f_1,\cdots ,f_n ) \widetilde = Mat_{2^n}(\R) \] and
     \[ \G_{n,n+1}:= \C(e_1, \cdots , e_n,f_1,\cdots ,f_n ) \widetilde = Mat_{2^n} (\C),   \] 
    where $e_i := a_i+b_i$ and $f_i := a_i-b_i$ for $i=1, \cdots, n$.  
    The geometric algebra $\G_{n,n+1}$ can be obtained as a real geometric algebra,
    \beq \G_{n,n+1}:= \R(e_1, \cdots, e_n,f_1, \cdots ,f_n,f_{n+1}),  \label{geospq} \eeq 
    where $f_{n+1}:= e_{1\cdots n}f_{1\cdots n}^\dagger i $.
 
 If $p+q \le 2n$, the $2^{2n}$ dimensional geometric algebras $\G_{p,q}$ can be considered as the real geometric algebras, 
  \[  \G_{p,q} := \R(e_1, \cdots , e_p, f_1, \cdots, f_q),     \]
   by noting that $(ie_i)^2=-1$ and $(i f_i)^2=1$. If $p>n$ or $q>n$, we simply change the required number of basis vectors $f_i$ or $e_i$, to $e_i$'s or $f_i$'s, respectively, by multiplying the $f_i$ or $e_i$ by $i$. 
     Additional $2^{2n+1}$-dimensional real geometric algebras can be obtained from
     (\ref{geospq}), by replacing an {\it even} number of the generating vectors by $i$ times that generator. For example,
     \[ \G_{n-2,n+3}:= \R(e_1, \cdots, e_{n-2}, i e_1,i e_2,f_1, \cdots, f_{n+1} )
     \widetilde = Mat_\C(2^n)  .    \]
     
     Below is a list of lower dimensional $2^{2n+1}$ geometric algebras that can be represented in terms of the matrix algebras     
     $Mat_2(\C)$, $Mat_{2^2}(\C)$, and $Mat_{2^3}(\C)$,  or vice-versa.
     
     \bigskip
     
    \centerline{$\bf  \G_{1,2}:= \C(e_1,f_1) \widetilde =Mat_2(\C) $} 
    
    \[ \G_{1,2} = \R(e_1,f_1,f_2), \quad \G_{3,0}= \R(e_1,i f_1, i f_2) ,  \]   
  
       \centerline{$\bf  \G_{2,3}:= \C(e_1,e_2,f_1,f_2) \widetilde =Mat_4(\C) $} 
       
       \[ \G_{2,3} = \R(e_1,e_2,f_1,f_2,f_3), \quad \G_{4,1}= \R(e_1,e_2,if_1, if_2,f_3) ,  \]
       \[ \G_{0,5} = \R(ie_1,ie_2,f_1,f_2,f_3) .  \]
       
        \centerline{$\bf  \G_{3,4}:= \C(e_1,e_2,e_3,f_1,f_2,f_3,f_4) \widetilde =Mat_8(\C) $} 
        
        \[ \G_{3,4} = \R(e_1,e_2,e_3,f_1,f_2,f_3,f_4), \quad \G_{5,2}= \R(e_1,e_2,e_3,if_1, if_2,f_3,f_4) ,  \]
        \[ \G_{7,0} = \R(e_1,e_2,e_3,i f_1,i f_2,i f_3,i f_4), \quad  \G_{1,6} = \R(e_1,i e_2,i e_3, f_1,f_2,f_3,f_4) .   \]
      A complete classification of geometric algebras is given in the
      paper, {\it Geometrization of the Real Number System} \cite{geoS2017}.


\begin{thebibliography}{}
  	\bibitem{TD1967} T. Dantzig, {\em Number: The Language of Science}, 4th edn. Free Press, New York 1967.
  	\bibitem{S08} G. Sobczyk, Geometric Matrix Algebra, 
  	{\em Linear Algebra and its Applications}, 429 (2008) 1163-1173.
  	 \bibitem{McKay} D.J.C. MacKay, {\it Good Error Correcting Codes based on Very Sparse Matrices}, IEEE Transactions on Information Theory, Vol: 45, Issue: 2,pp. 399-431, Mar 1999. 
  	 \bibitem{kao2010} D. Kao, {\it Representations of the Symmetric Group},
  	 http://www.math.uchicago.edu/~may/VIGRE/VIGRE2010/REUPapers/Kao.pdf
  	http://www.inference.org.uk/mackay/mncN.pdf
  	 \bibitem{JL2001} G. James, M. Liebeck {\it Representations and Characters of Groups}, Second Edition 2nd Edition, Cambridge 2001.  
  	  \bibitem{FH1991}W. Fulton, J. Harris, {\it Representation Theory:
  	 	A First Course}, Springer-Verlag 1991.
  	  	\bibitem{Weyl1950} H. Weyl, {\it Theory of Groups and Quantum Mechanics}, Dover Publications, Inc. 1950.
  	\bibitem{hypre2017} G. Sobczyk, {Hyperbolic Numbers Revisited}, 
  	  	http://www.garretstar.com/hyprevisited12-17-2017.pdf
   \bibitem{SNF} G. Sobczyk, {\em New Foundations in Mathematics: The Geometric Concept of Number},
  	  	\newblock Birkh\"auser, New York 2013.		 \begin{verbatim} http://www.garretstar.com/ \end{verbatim}
   \bibitem{Kun2015} Kunle Adegoke, Olawanle Layeni, Rauf Giwa, Gbenga Olunloyo, {\it The Standard Representation of the Symmetric Group
  	  		 	$S_n$ over the Ring of Integers}, Turkish Journal of Analysis and Number Theory, 2015, Vol. 3, No. 5, 126-127.
  	  		 Available online at http://pubs.sciepub.com/tjant/3/5/3 
   \bibitem{newman67} M. Newman, {\it Two Classical Theorems on Commuting Matrices}, JOURNAL OF RESEARCH of the National Bureau of Standards-B, Mathematics and Mathematic al Physics, Vol. 71 B, Nos. 2 and 3, April-September 1967. https://nvlpubs.nist.gov/nistpubs/jres/71b/jresv71bn2-3p69\_a1b.pdf		 
  	 \bibitem{geoS2017} G. Sobczyk, {\it Geometrization of the Real Number System}, July, 2017.
  	 \begin{verbatim}http://www.garretstar.com/geonum2017.pdf\end{verbatim}	
  	
 	\bibitem{S2} G. Sobczyk, {\it The Generalized Spectral Decomposition
 		of a Linear Operator}, The College Mathematics
 	Journal, 28:1 (1997) 27-38.
  	\bibitem{S0} G. Sobczyk, {\it The missing spectral basis in algebra and number theory}, The American Mathematical Monthly 108 April 2001, pp. 336-346.
  
  	 \end{thebibliography}
  	\end{document}